\documentclass[a4paper,10pt,reqno]{amsart}
\usepackage{verbatim}
\usepackage{amssymb}
\usepackage{mathrsfs}

\usepackage{enumerate}
\usepackage[active]{srcltx}
\numberwithin{equation}{section}
\usepackage{tikz}
\usetikzlibrary{backgrounds}
\usetikzlibrary{positioning}
\usetikzlibrary{snakes}
\usetikzlibrary {arrows.meta}

 \makeatletter
  \@namedef{subjclassname@2020}{%
  \textup{2020} Mathematics Subject Classification}
   \makeatother

\usepackage{t1enc}
\usepackage[utf8x]{inputenc}

\newtheorem{theorem}{Theorem}[section]
\newtheorem{lemma}[theorem]{Lemma}

\newtheorem{fact}[theorem]{Fact}
\newtheorem{proposition}[theorem]{Proposition}

\newtheorem{corollary}[theorem]{Corollary}
\newtheorem{claim}{Claim}[theorem]

\theoremstyle{definition}
\newtheorem{definition}[theorem]{Definition}
\theoremstyle{remark}
  
\newcommand{\mc}[1]{\mathcal{#1}}

\newcommand{\mf}[1]{\mathfrak{#1}}
\newcommand{\setm}{\setminus}
\newcommand{\empt}{\emptyset}
\newcommand{\subs}{\subset}
\newcommand{\llambda}{{{\omega}_1}}

\newcommand{\mysubsection}[1]{\smallskip\noindent{\bf #1}\smallskip}

\DeclareMathOperator{\Fn}{Fn}
\DeclareMathOperator{\supp}{supp}

\DeclareMathOperator{\dom}{dom}

\DeclareMathOperator{\weight}{w}

\DeclareMathOperator{\PNA}{PNA}

\newcommand{\<}{\left\langle}
\renewcommand{\>}{\right\rangle}

\newcommand{\stickT}{%
\setbox255=\hbox{\raise1ex\hbox{$\hspace{0.2pt}\,\bullet\,$}}
\mathord{\rlap{\hbox to\wd255{\hss\hbox{$|$}\hss}}
\box255}
}
\newcommand{\stickS}{%
\setbox255=\hbox{\raise0.6ex\hbox{$\scriptstyle\bullet$}}
\mathord{\rlap{\hbox to\wd255{\hss\hbox{$\scriptstyle|$}\hss}}
\box255}
}
\newcommand{\stick}{{\mathchoice{\stickT}{\stickT}{\stickS}{\stickS}}}
\DeclareMathOperator{\nw}{nw}
\newcommand{\superstick}{\stick_{\mbox{\tiny SUPER}}}
\DeclareMathOperator{\den}{d}
\DeclareMathOperator{\hden}{z}
\DeclareMathOperator{\hL}{h}

\author[I. Juh{\'a}sz]{Istv{\'a}n Juh{\'a}sz}
\thanks
  {
   }
\address
      { Alfr{\'e}d R{\'e}nyi Institute of Mathematics,
E{\"o}tv{\"o}s Lor{\'a}nd Research Network
}
\email{juhasz@renyi.hu}

\author[L. Soukup]{Lajos Soukup}
\thanks
  {
   }
\address
      { Alfr{\'e}d R{\'e}nyi Institute of Mathematics,
      E{o}tv{o}s Lor{\'a}nd Research Network
}
\email{soukup@renyi.hu}

\author[Z. Szentmikl\'ossy]{Zolt\'an Szentmikl\'ossy}
\address{E{o}tv{o}s University of Budapest}
\email{szentmiklossyz@gmail.com}


\title{The class $C({\omega}_1)$ and countable net weight}
\thanks{The research on and preparation of this paper was
supported by  NKFIH grant   K129211}
\date{\today}

\begin{document}

\begin{abstract}
Hart and Kunen in \cite{HKsup} and, independently, Ríos-Herrejón in \cite{rios2023weakly} defined and studied
the class $C({\omega}_1)$ of topological spaces $X$ having the property that for every neighborhood
assignment $\{U(y) : y \in Y\}$ with $Y \in [X]^{\omega_1}$ there is $Z \in [Y]^{\omega_1}$ such that $$Z \subs \bigcap \{U(z) : z \in Z\}.$$
It is obvious that spaces of countable net weight, i.e. having a countable network, belong to this class.
In this paper we present several independence results concerning the relationships of these two and several other natural classes
that are sandwiched between them, thus
clarifying some of the main problems that were raised  in \cite{rios2023weakly}.

In particular, we prove that the continuum hypothesis, in fact a weaker combinatorial principle called {\em super stick},
implies that every regular space in $C({\omega}_1)$ has countable net weight, answering a question that was raised by Hart and Kunen in both
\cite{HKsup} and \cite{HKsupdup}.
\end{abstract}

\maketitle

\section{Introduction}

In this paper the notation and terminology concerning set theory follows \cite{Kunen} and concerning
cardinal functions follows \cite{Juh10}.

The concept of weakly separated spaces, a common weakening of both left- and right-separated ones,
was introduced by Tkachenko in \cite{Tka78}. We recall that $X$ is {\em weakly separated} iff there is
a neighborhood assignment $U : X \to \tau(X)$ on $X$ such that if $x, y$ are distinct points in $X$
then $x \notin U(y)$ or $y \notin U(x)$. Also, he defined $R(X)$ as the supremum of the cardinalities
of all weakly separated subspace of $X$. The related "hat version"
$\widehat{R}(X)$ is the minimal cardinal ${\mu}$
such that $X$ does not contain a weakly separated subspace of size ${\mu}$.

Clearly, if $X$ has a countable network, i.e. $\nw(X)={\omega}$,
then $R(X)={\omega}$ (and even $R(X^{\omega})={\omega}$), so the natural question was raised in \cite{Tka78}
if the converse of these implications hold.
Several consistent counterexamples to these have been provided since then, see e.g. in timeline:
\cite{HaJu79,Cie83,Tod84PartProb,JuSoSz94WhatMakes}, and \cite{HKsup},
but no ZFC example of a space $X$ with $R(X) = \omega < nw(X)$ is currently known.

More recently the class  $C({\omega}_1)$ of spaces, that is sandwiched between those of countable net weight
and those satisfying $R(X^{\omega})={\omega}$,
was introduced in \cite{HKsup} and later, apparently independently, in \cite{rios2023weakly}.
Actually, in \cite{HKsup} the spaces in this class are called {\em super hereditarily good}, in short: {\em suHG}.
(It was only after the first version of our present paper appeared on the arXiv
that \cite{HKsup} came to our attention.)

In \cite[Question 4.21]{rios2023weakly}, the author
posed the following related general question:
how does $X \in C({\omega}_1)$ relate to the assumptions
$R(X^{\omega})={\omega}$ and $\nw(X)={\omega}$? (Actually,
$R(X^{\omega})$ is denoted by $R^*(X)$ in \cite{rios2023weakly}.)

Before giving the definition of $C({\omega}_1)$, or more generally of $C({\kappa})$
for any uncountable cardinal $\kappa$, we introduce some terminology and notation.
Given any topological space $X$, a {\em partial neighborhood assignment on $X$} is a function $U$
whose domain is a subset $Y$ of $X$ and $y \in U(y) \in \tau(X)$ for any $y \in Y$.
We shall denote by $\PNA(X, \kappa)$ the collection of all partial neighborhood assignments $U$
on $X$ for which $|\dom(U)| = \kappa$.

The following general definition was introduced in \cite{rios2023weakly}.

\begin{definition}\label{df:Cakappa}
The topological space $X$ has property $C({\kappa})$, denoted as  $X\in C({\kappa})$, iff
for
any $U \in \PNA(X, \kappa)$
there exists a set $Y\in {[\dom(U)]}^{{\kappa}}$ such that
\begin{displaymath}
Y\subs\bigcap_{y\in Y}U(y).
\end{displaymath}
\end{definition}

We note that in the paper \cite{HKsup}, predating \cite{rios2023weakly},
only the class suHG that is identical with $C({\omega}_1)$ was considered.

In
\cite[Proposition 3.14, Corollary 4.20, Question 4.21]{rios2023weakly} the author
made observations (1) and (2) below and raised question (3) that prompted our present work.

 {\em Assume that  $X$ is a topological space and
${\kappa}$ is a regular cardinal.
\begin{enumerate}[(1)]
\item If $\nw(X)<{\kappa}$, then $X\in C({\kappa})$.
\item If ${\omega}<{\kappa}<\widehat R(X^{\omega})$, then
$X \notin C({\kappa})$.
\end{enumerate} }
\begin{enumerate}[(1)]
\item[\it (3)] {\em When does $\widehat R(X^{\omega})\le {\kappa}\le \nw(X)$
imply $X\in C({\kappa})$?}
\end{enumerate}

We next introduce some additional properties that
 allow  for a more detailed  analysis of our subject.

\begin{definition}\label{df:calK}
Given any partial neighborhood assignment
$U$ on a space $X$, we let  ${\tau}_U$ denote  the topology
on $\dom(U)$ generated by the family $$\{U(x) \cap \dom(U) : x \in \dom(U)\}.$$
The resulting space $\<\dom(U),{\tau}_U\>$ is denoted by $X_U$.
Clearly, $\weight(X_U) \le |\dom(U)|$.
\end{definition}

\begin{definition}\label{df:Kakappa}
The space $X$ has property $K({\kappa})$ iff for
any $U \in \PNA(X, \kappa)$
we have $\nw(X_U)<{\kappa}$.
   \end{definition}

   \begin{definition}\label{df:Nakappa}
      $X$ has property $N({\kappa})$ iff
      $\nw(Y)<{\kappa}$ for each subspace $Y\in {[X]}^{{\kappa}}$.
      \end{definition}

Clearly,
\begin{displaymath}
\nw(X)<{\kappa} \to\  X\in N({\kappa})\ \to\ X\in K({\kappa})
\stackrel{cf({\kappa})={\kappa}}{\xrightarrow{\hspace*{1cm}}}\ X\in C({\kappa})\to
\widehat{R}(X^{\omega})\le{\kappa}.
\end{displaymath}
In this paper on one hand we find conditions ensuring that some of these implications are consistently reversible,
and on the other we show that these results are not provable in ZFC.

In order to orient the reader, we first list the related results known to us.

\begin{enumerate}[(a)]
\item \label{HaJu-DiMa} Hajnal and Juhász (\cite[Theorem 3]{HaJu79}):  If both $\diamondsuit_{{\omega}_2}^{\omega}$ and $MA({\omega}_1)$ hold, then
there exists a space $X\in N({\omega}_1) \cap T_2$ with $|X|={\omega}_2$ and $\nw(X)>{\omega}$.
\item \label{Ci-fo} Ciesielski (\cite{Cie83}):
In some CCC generic extension,  $\mathfrak{c} = 2^{\omega}>{\omega}_1$ and there exists a 0-dimensional
space $X\in N({\omega}_1) \cap T_2$
with $|X|=\weight(X)=\nw(X)={\omega}_2$.
\item \label{boo}
%
Todorcevic (\cite[Theorem 3.5]{Tod84PartProb} and the first
sentence of its proof): If $\mf b={\omega}_1$, then
there is a 0-dimensional $T_2$-space $X$ with $R(X^\omega)={\omega}$ such that
$|X|=\weight(X)={\omega}_1$ and there is a neighborhood assignment $U:X\to {\tau}(X)$ such that   for each $A\in {[X]}^{{\omega}_1}$
there is $\{x,y\}\in [A]^2$ with $x\notin U(y)$ and $y\notin U(x)$
(and so  $X\notin C({\omega}_1)$).

\item \label{blarge}
Todorcevic (\cite[Theorem 3.7]{Tod84PartProb}): If $\mf b>{\omega}_2$ then
there is 0-dimensional $T_2$-space $X \in N({\omega}_1)$ with  $|X|={\omega}_2$ and
$\nw(X)>{\omega}$.
\item  Juhász et al.   \cite[Theorem 3.5]{JuSoSz94WhatMakes}:  \label{jusosz-gen}
In a certain CCC generic extension,  there is a 0-dimensional $T_2$-space $X$
with $R(X^{\omega})={\omega}$ but $X\notin C({\omega}_1)$.
\item Juhász et al. \cite[Corollary 3.5]{JuSoSz96countablenetwork}\label{juhoszo2}:
$MA_{\omega_1}$
implies that if $X$ is any $T_2$-space with  $|X|+\weight(X)\le {\omega}_1$
then  $\nw(X)={\omega}$ iff $R(X^{\omega})={\omega}$.

\item \label{HKtm1.4}
Hart and Kunen (\cite[Theorem 1.4]{HKsup}): It is consistent with ZFC to have MA
and a first countable space $X\in C({\omega}_1) \cap T_3$ such that
$$|X| = \weight(X) = \nw(X) = \mathfrak{c} = \omega_2.$$

\end{enumerate}

\smallskip

In this paper we shall prove the following relevant new results:

\begin{theorem}\label{tm:ch}
The combinatorial principle $\superstick$, an obvious consequence of CH, implies
$\nw(X)={\omega}$ for any $X\in C({\omega}_1) \cap T_3$.
\end{theorem}

As we mentioned in the abstract, the question if CH implies that all spaces in $C({\omega}_1) \cap T_3$
have countable net weight was raised in both  \cite{HKsup} and \cite{HKsupdup}.

\begin{theorem}\label{tm:cohen}
   If $V \models CH$ and $\kappa$ is any cardinal, then after adding $\kappa$ Cohen reals we get
   $C({\omega}_1) \cap T_3 \subset K({\omega}_1)$ in the extension, i.e.

   \begin{displaymath}
   V^{{\Fn}({\kappa},2)}\models
   C({\omega}_1) \cap T_3 \subset K({\omega}_1).
   \end{displaymath}
   \end{theorem}

   \begin{theorem}\label{tm:example}
      It is consistent that $2^{\omega}$ is arbitrarily large and
      there is a 0-dimensional $T_2$-space $X\in K({\omega}_1)$
      such that  $|X|=\nw(X)={\omega}_1$ and $\weight(X)=2^{\omega}$,
      So, in this case $X$ is a witness for
      $$K({\omega}_1) \cap T_3 \not\subset N({\omega}_1).$$
      \end{theorem}

      \begin{theorem}\label{tm:stac}
         It is consistent that there is a 0-dimensional $T_2$-space $X=\<{\omega}_1,{\tau}\>$ such that
         $\weight(X)={\omega}_1$ and
         \begin{displaymath}
         \{Y\subs X: \nw(Y)={\omega}\}=NS({\omega}_1),
         \end{displaymath}
         where $NS(\omega_1)$ is the ideal of non-stationary subsets of $\omega_1$.
         An easy consequence is that this $X\in C({\omega}_1)\setm K({\omega}_1)$.
         \end{theorem}

\begin{figure}[h]\label{masodikabra}
   \begin{tikzpicture}[scale=0.85]
      \node[draw] (a) at (0,0) {\small $\nw(X)={\omega}$};
      \node[draw] (b) at (3,0) {\small $X\in N({\omega}_1)$};
      \node[draw] (c) at (6,0) {\small $X\in K({\omega}_1)$};
      \node[draw] (d) at (9,0) {\small $X\in C({\omega}_1)$};
      \node[draw,dotted] (e) at (12,0) {\small $R(X^{\omega})={\omega}$};
      \path[draw,-{Stealth[length=2mm]}]  (a) -- (b) (b) edge (c) (c) edge    (d)  (d) edge (e);
     \path (d) edge[out=145,in=75,-{Stealth[length=2mm]},dashed]  node[left] {\small Thm \ref{tm:cohen} }  (c) ;
     \path (e) edge[out=110,in=20,-{Stealth[length=2mm]},dashed]  node[above] {\small  \eqref{juhoszo2} }  (c) ;
   \path[draw,dotted] (1.5,1) -- (1.5,-1) node[below] {\small \eqref{Ci-fo},\eqref{blarge}};
   \path[draw,dotted] (4.5,1) -- (4.5,-1) node[below] {\small Thm \ref{tm:example}};
   \path[draw,dotted] (7.5,1) -- (7.5,-1) node[below] {\small Thm \ref{tm:stac}};
   \path[draw,dotted] (10.5,1) -- (10.5,-1) node[below] {\small \eqref{boo},\eqref{jusosz-gen}};
   \path[draw,-{Stealth[length=2mm]},dashed] (d) -- (9,1.8) -- node[above] {\small Thm \ref{tm:ch}}  (0,1.8) --(a);
   \end{tikzpicture}
   \caption{Results concerning  regular spaces}
      \end{figure}
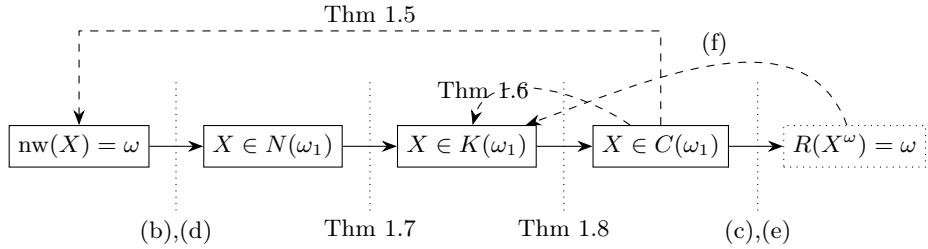

Figure \ref{masodikabra} provides a summary of our knowledge regarding
regular spaces when ${\kappa}={\omega}_1$.

We are very interested in the following problems that remained open.
\begin{enumerate}[(1)]
\item Does CH imply that every  $X\in C({\omega}_1) \cap T_2$ has countable net weight?
\item Is it true that if $V \models CH$ then, for every cardinal $\kappa$,
\begin{displaymath}
    V^{{\Fn}({\kappa},2)}\models
   C({\omega}_1) \cap T_3 \subset N({\omega}_1)?
   \end{displaymath}

\item Is it provable in ZFC that every uncountable space $X\in C({\omega}_1) \cap T_3$ has
an uncountable subspace of countable net weight?
\end{enumerate}

We also mention that in \cite{HKsup} the authors asked if PFA implies the existence of a space
$X \in C({\omega}_1) \cap T_3$ with $\nw(X) > \omega$; this question remains open as well.

\section{Conditions implying $\nw(X) < \kappa$ for $X \in C(\kappa)$}

We start by recalling that $\widehat{z}(X)$ (resp. $\widehat{h}(X)$) is the minimal cardinal ${\mu}$
such that the space $X$ does not contain a left-separated (resp. right-separated) subspace of size ${\mu}$.

\begin{theorem}\label{tm:kappaomega}
Assume that $X\in C({\kappa}) \cap T_3$ and  ${\kappa}^{<{\mu}}={\kappa}$, where
$\mu = \widehat \hden(X)$. Then $\nw(X)<{\kappa}$.
In particular,  any hereditarily separable $T_3$-space  $X\in C(\mathfrak{c})$
has net weight less than $\mathfrak{c}$.
\end{theorem}

\begin{proof}
Assume that the theorem fails and $X \in C({\kappa}) \cap T_3$ is a counterexample, i.e. $\nw(X)\ge {\kappa}$.
By ${\kappa}^{<{\mu}}={\kappa}$ we may fix an enumeration
 $\{S_{\zeta}:{\zeta}<{\kappa}\}$ of ${[{\kappa}]}^{<{\mu}}$.
We then define $U = \{\<x_{\alpha},U_{\alpha}\>:{\alpha}<{\kappa}\} \in \PNA(X, \kappa)$
by transfinite recursion on ${\alpha}< \kappa$ as follows.

Assuming that ${\alpha}< \kappa$ and $\<x_\beta, U_\beta\>$ have been defined for $\beta < \alpha$, we let
\begin{displaymath}
\mc N_{\alpha}=
\big\{\overline {\{x_{\nu}:{\nu}\in S_{\zeta}\}}:{\zeta}<{\alpha} \text{ and } S_{\zeta}\subs {\alpha}\}
\cup \{\{x_{\nu}\}:{\nu}<{\alpha}\big\}.
\end{displaymath}
Now, $|\mc N_{\alpha}| < \kappa$ implies that $\mc N_{\alpha}$ is not a network for $X$. So it follows from  $X \in T_3$ that
there exist $x_{\alpha}\in X$ and $U_{\alpha}\in \tau(X)$ with
$x_{\alpha}\in U_{\alpha}$
such that  $x_{\alpha}\in T$ implies $T \not\subs  \overline {U_{\alpha}}$ for all $T\in \mc N_{\alpha}$.

Now, $X\in C({\kappa})$ implies that there  is $I\in {[{\kappa}]}^{{\kappa}} $ such that
\begin{displaymath}\tag{$*$}
\{x_{\nu}:{\nu}\in I\}\subs \bigcap\{U_{{\nu}}:{\nu}\in I\}.
\end{displaymath}
But we have $\den(\{x_{\alpha}:{\alpha}\in I\})< \mu = \widehat \hden(X)$,
hence there is ${\xi}<{\kappa}$ such that $S_{\xi}\subs I$ and
$\{x_{\alpha}:{\alpha}\in S_{\xi}\}$ is dense in $\{x_{\alpha}:{\alpha}\in I\}$.
We may then pick ${\alpha}\in I$ with ${\xi}<{\alpha}$ such that $S_{\xi}\subs {\alpha}$.
Indeed, this is because ${\kappa}^{<{\mu}}={\kappa}$ implies $cf(\kappa) \ge \mu$,
hence $S_{\xi}$ is bounded in $I$.

But then, on one hand  $x_{\alpha}\in \overline{\{x_{\nu}:{\nu}\in S_{\xi}\}} \in \mc N_{\alpha}$, and so
$\overline{\{x_{\nu}:{\nu}\in S_{\xi}\}}\not\subs \overline{U_{\alpha}},$
on the other hand $\{x_{\nu}:{\nu}\in S_{\xi}\} \subs U_{\alpha}$ by $(*)$, so
$\overline{\{x_{\nu}:{\nu}\in S_{\xi}\}}\subs \overline{U_{\alpha}}$.
This is a contradiction that completes the proof of the theorem.
\end{proof}

Since $\widehat{\hL}(X)$ is in a certain sense dual to $\widehat{\hden}(X)$, our next result may be
considered dual to the previous one.

\begin{theorem}\label{tm:weightkappalambda}
Assume that $X\in C({\kappa}) \cap T_3$ and  ${\kappa}^{<{\mu}}={\kappa} = \weight(X)$, where
$\mu = \widehat \hL(X)$. Then $\nw(X)<{\kappa}$.
In particular, for any hereditarily Lindelöf $T_3$-space $\,X\in C(\mathfrak{c})$ of weight $\mathfrak{c}$ we have
$\nw(X) < \mathfrak{c}$.
  \end{theorem}

  \begin{proof}[Proof of Theorem \ref{tm:weightkappalambda}]
Assume that the theorem fails and $X$ is a counterexample.
Now, $\mu = \widehat h(X)$ means that for any $\mathcal{U} \subs \tau(X)$
there is $\mathcal{V} \in [\mathcal{U}]^{< \mu}$ such that $\cup \mathcal{U} = \cup \mathcal{V}$.
Consequently, $\weight(X) = \kappa = {\kappa}^{<{\mu}}$ imply $|\tau(X)| = \kappa$ as well.
Also, if $\mc H$ denotes the family of all $G_{< \mu}$-sets in $X$, then $|\mc H|= |{\tau}(X)|^{<{\mu}}\ = {\kappa}$ as well.

Let $\{H_{\zeta}:{\zeta}<{\kappa}\}$ be an enumeration of $\mc H$.
We will again define a sequence $\{\<x_{\alpha},U_{\alpha}\>:{\alpha}<{\kappa}\}$
by transfinite recursion on ${\alpha}< \kappa$.

Assuming that ${\alpha}< \kappa$ and $\{\<x_\zeta, U_\zeta\> : \zeta < \alpha \}$
has been defined, we let
\begin{displaymath}
\mc N_{\alpha}=
\big\{H_{\zeta} : {\zeta}<{\alpha}\big\} \cup \big\{\{x_{\zeta}\}:{\zeta}<{\alpha}\big\}.
\end{displaymath}
Then $|\mc N_{\alpha}| < \kappa$, so it follows from $\nw(X) = {\kappa}$ and $X \in T_3$ that
there exist $x_{\alpha} \in U_{\alpha}\in \tau(X)$
such that for all $\zeta < \alpha$ we have $x_ \alpha \ne x_\zeta$, moreover
$H_\zeta \not\subs  \overline {U_{\alpha}}$ if $x_{\alpha} \in H_\zeta$.

By $X\in C({\kappa})$ there  is $I\in {[{\kappa}]}^{{\kappa}} $ such that
\begin{displaymath}\tag{$\star$}
\{x_{\nu}:{\nu}\in I\}\subs \bigcap\{U_{{\nu}}:{\nu}\in I\}\subs \bigcap\{\overline{U_{{\nu}}}:{\nu}\in I\}.
\end{displaymath}
Since $\mu = \widehat h(X)$, there is $J\in {[I]}^{<{\mu}}$ for which
\begin{displaymath}
\bigcup_{{\nu}\in I} (X\setm \overline {U_{\nu}})=
\bigcup_{{\nu}\in J} (X\setm \overline {U_{\nu}}),
\end{displaymath}
or equivalently,
\begin{displaymath}
\bigcap_{{\nu}\in I} \overline {U_{\nu}} =
\bigcap_{{\nu}\in J} \overline {U_{\nu}}.
\end{displaymath}
Now, $\bigcap_{{\nu}\in J} U_{\nu} \in \mc H$, hence $\bigcap_{{\nu}\in J} U_{\nu} = H_\xi$
for some ${\xi}<{\kappa}$. Let us pick $\alpha \in I$ with $\xi < \alpha$. Then on one hand
$$x_\alpha \in \bigcap_{{\nu}\in I} U_{\nu} \subs \bigcap_{{\nu}\in J}
 U_{\nu} = H_\xi$$
implies $H_\xi \not\subs  \overline {U_{\alpha}}$. But on the other hand we have
\begin{displaymath}
   H_\xi = \bigcap_{{\nu}\in J} {U_{\nu}}\subs \bigcap_{{\nu}\in J} \overline {U_{\nu}} = \bigcap_{{\nu}\in I} \overline {U_{\nu}} \subs \overline {U_{\alpha}}.
\end{displaymath}
This contradiction completes the proof.
\end{proof}
We do not know if the assumption $\weight(X) = \kappa$ can be dropped from this result.

Clearly, $X\in C({\kappa})$ implies both $\widehat z(X) \le \kappa$ and $\widehat h(X) \le \kappa$.
Moreover, if $X$ is $T_3$ then $d(X) < \widehat z(X) \le \kappa$ and  $\weight(X) \le 2^{d(X)} \le \kappa$
if $\kappa^{< \kappa} = \kappa$. Thus, both theorems above imply the following.

\begin{corollary}\label{co:kík}
If $\kappa^{< \kappa} = \kappa$ then $\nw(X) < \kappa$ for every $X\in C({\kappa}) \cap T_3$.
\end{corollary}

Our next result is a strengthening of Corollary \ref{co:kík} in that it yields the same conclusion
from a combinatorial principle that is weaker than $\kappa^{< \kappa} = \kappa$.

\begin{definition}\label{df:superstick}
$\Sigma(\kappa)$ is the statement that there is a family $\mc S\subs [\kappa]^{< \kappa}$
with $|\mc S|= \kappa$ such that
\begin{displaymath}
\forall I \in [\kappa]^\kappa\, \exists J\in [I]^\kappa\
\big(\forall  T \in {[J]}^{< \kappa}\ \exists S\in \mc S \ (T \subs S \subs I)\big).
\end{displaymath}
\end{definition}

Clearly, $\kappa^{< \kappa} = \kappa$ implies $\Sigma(\kappa)$ because then we may take $\mc S = [\kappa]^{< \kappa}$.
In particular, CH implies $\Sigma(\omega_1)$. Actually,
$\Sigma(\omega_1)$ was introduced by Primavesi in
\cite{primavesi2011guessing}, where it was called {\em super stick}, denoted as  $\superstick$.

In \cite[Theorem 3.1]{Chen17} it is claimed that $\neg CH + \superstick$ is consistent modulo the existence
of an inaccessible cardinal.
Garti and Shelah  in \cite{gartishelah2023tiltan} also claim the consistency of $\neg CH + \superstick$.

Note that $\kappa^{< \kappa} = \kappa$ also implies that $\kappa$ is a regular cardinal.
\begin{theorem}\label{tm:supeclub}

If $\kappa$ is a regular cardinal and $\Sigma(\kappa)$ holds then $\nw(X) < \kappa$ for every
$X\in C({\kappa}) \cap T_3$.
\end{theorem}

\begin{proof}[Proof of Theorem \ref{tm:supeclub}]
Our proof is again indirect, so assume that $X\in C({\kappa}) \cap T_3$ with $\nw(X) \ge \kappa$.

Let $\mc S=\{S_{\zeta}:{\zeta}< \kappa \}$ be a family that witnesses $\Sigma(\kappa)$.
We shall define a sequence $\{\<x_{\alpha},U_{\alpha}\>:{\alpha}<{\kappa}\}$
such that $U = \{\<x_{\alpha},U_{\alpha}\>:{\alpha}<{\kappa}\} \in \PNA(X, \kappa)$
by transfinite recursion on ${\alpha}<{\kappa}$.

Having defined $\{\<x_\zeta , U_\zeta\> : \zeta < \alpha\}$
 we let
\begin{displaymath}
\mc N_{\alpha}=
\big\{\overline {\{x_{\nu}:{\nu}\in S_{\zeta}\}} : {\zeta}<{\alpha} \text{ and } S_{\zeta}\subs {\alpha}\}
\cup \{\{x_{\zeta}\}:{\zeta}<{\alpha}\big\}.
\end{displaymath}
Since $\nw(X) \ge \kappa$ and $X$ is $T_3$,
we may choose $x_{\alpha}\in X$ and
$U_\alpha \in {\tau}(X)$ with $x_{\alpha}\in U_{\alpha}$
such that  $x_{\alpha}\in T$ implies $T \not\subs  \overline {U_{\alpha}}$ for all $T\in \mc N_{\alpha}$.

By $X\in C(\kappa)$ there  is $I\in [\kappa]^\kappa$ such that
\begin{displaymath}\tag{$*$}
\{x_{\nu}:{\nu}\in I\}\subs \bigcap\{U_{{\nu}}:{\nu}\in I\},
\end{displaymath}
and $\mc S$  being a $\Sigma(\kappa)$-family,
there is $J\in {[I]}^{\kappa}$ such that
for every $T \in {[J]}^{< \kappa}$ there is ${\zeta} < \kappa$
such that $T \subs S_{\zeta} \subs I$.

But $X\in C({\kappa})$ also implies that $\den(\{x_\nu : \nu \in J\}) < \kappa$,
hence there is $T \in {[J]}^{< \kappa}$ for which $\{x_\nu : \nu \in T\}$ is dense in $\{x_\nu : \nu \in J\}$.
Then $T \subs S_{\xi}\subs I$ for some ${\xi}<\kappa$. Since $\kappa$ is regular, there is
${\alpha}\in J$ such that $\xi < \alpha$ and $S_{\xi} \subs {\alpha}$, hence $\overline {\{x_{\nu}:{\nu}\in S_{\xi}\}} \in \mc N_{\alpha}$.
Then we would have both (1) and (2) below:
\begin{enumerate}[(1)]
\item $x_{\alpha}\in \overline{\{x_{\nu}:{\nu}\in S_{\xi}\}} \in \mc N_{\alpha}$, and so
$\overline{\{x_{\nu}:{\nu}\in S_{\xi}\}}\not\subs \overline{U_{\alpha}}. $

\smallskip

\item $\{x_{\nu}:{\nu}\in S_{\xi}\} \subs U_{\alpha}$ by $(*)$, so
$\overline{\{x_{\nu}:{\nu}\in S_{\xi}\}}\subs \overline{U_{\alpha}}$.
\end{enumerate}
This contradiction completes the proof.
\end{proof}

Since the methods used in the proofs of the above three theorems look very similar,
it is natural to wonder if they could be deduced from a single result as corollaries.
In fact, this question was raised by the referee but, alas, we could not provide
such a single source for our three theorems..

\section{Cohen forcing and $C({\omega}_1) \cap T_3 \subs K(\omega_1)$}

The aim of this section is to prove Theorem \ref{tm:cohen} from the introduction, i.e.
that if $V \models CH$ then $V^{{\Fn}({\kappa},2)}\models  C({\omega}_1) \cap T_3 \subs K(\omega_1)$
for any cardinal $\kappa$.
Since $U \in \PNA(X, \omega_1)$ implies $\weight(X_U) \le \omega_1$, the theorem clearly follows if we can
show that, in $V^{{\Fn}({\kappa},2)}$, we have $\nw(X) = \omega$ for any $X \in C({\omega}_1) \cap T_3$
such that $|X| = \weight(X) = \omega_1$.

Now, if our ground model $V \models CH$ then $V^{{\Fn}({\omega_1},2)}\models  CH$ as well, hence  our first step
is to show that if $X \in V$ with $|X| = \weight(X) = \omega_1$ then $V^{{\Fn}({\kappa},2)}\models X \in C({\omega}_1)$ implies
$V^{{\Fn}({\omega_1},2)}\models X \in  C(\omega_1)$ as well.
This will actually follow from an analogous result on graphs that might be of independent interest.
A graph is a pair $\<S, E\>$, where $E \subs [S]^2$. (We use $S$ to denote the set of vertices
instead of the customary $V$ because $V$ now denotes our ground model.)

\begin{definition}\label{df:gra-cw1}
The graph $\<S, E\>$ with $|S| = \kappa$ has property  $\Gamma(\kappa)$ iff every spanned subgraph of it on $\kappa$ vertices
contains a clique of size $\kappa$, i.e.
\begin{displaymath}
  \forall A\in {[S]}^{\kappa}\, \exists \, B \in {[A]}^{\kappa}\, \big( [B]^2 \subs E \big).
\end{displaymath}
\end{definition}

\begin{theorem}\label{tm:cohen-graph}
   If $\<S, E\>$ with $|S| = \omega_1$ is a graph in our ground model $V$ then, for every cardinal $\kappa > \omega_1$,
   $$V^{{\Fn}({\kappa},2)}\models \<S, E\> \in \Gamma({\omega}_1)$$
implies
   $$V^{{\Fn}({\omega}_1,2)}\models \<S, E\> \in \Gamma({\omega}_1).$$
   \end{theorem}

   \begin{proof}[Proof of Theorem \ref{tm:cohen-graph}]
      Clearly, it suffices to show this for $S = \omega_1$. So, assume that
      \begin{equation}\label{eq:1}\notag
      1\Vdash_{{\Fn}({\kappa},2)} \<\omega_1, E\> \in \Gamma({\omega}_1),
      \end{equation}
To prove the theorem it is enough to show that
      for any $p \in {\Fn}({\omega}_1,2)$, if
$$p\Vdash_{{\Fn}({\omega}_1,2)} \dot A = \{\dot x_{\alpha}:{\alpha}<{\omega}_1\}\in {[S]}^{{\omega}_1},$$
then there is $r \in {\Fn}({\omega}_1,2)$ with $r \le p$ such that
$$r \Vdash_{{\Fn}({\omega}_1,2)}  \,\exists B \in [\dot A]^{\omega_1}\,([B]^2 \subs E).$$

To find $r$, we first fix
for each ${\alpha}<{\omega}_1$  a condition  $p_{\alpha}\in {\Fn}({\omega}_1,2)$ with $p_{\alpha}\le p$
      such that $p_{\alpha}$ decides the value of $\dot x_{\alpha}$.
       i.e.
      \begin{displaymath}
      p_{\alpha}\Vdash \dot x_{\alpha}= y_{\alpha}      \end{displaymath}
      for some $y_{\alpha}\in \omega_1$.
We can find  $I\in {[{\omega}_1]}^{{\omega}_1}$ such that
      $\{\dom(p_{\alpha}):{\alpha}\in I\}$ forms a $\Delta$-system with some kernel $D$ and
      $p_{\alpha}\restriction D = p_{\beta}\restriction D$ for all ${\alpha}, \beta \in I$.

Then we may pick $J\in {[I]}^{{\omega}_1}$ such that for
$K = \bigcup\{\dom(p_{\alpha}):{\alpha}\in J\}$ we have $|K| = |\omega_1 \setm K| = {\omega}_1$.
Clearly, $p \in {\Fn}(K, 2)$.

Now, consider the following ${\Fn}(K,2)$-name:
\begin{displaymath}
      \dot L=\{\<p_{\alpha}, {\alpha}\>:{\alpha}\in J\}.
      \end{displaymath}
      Then for any $M\subs {\kappa}$ with $K \subs M$ we have
      \begin{equation}\label{eq:pforce}
      p\Vdash_{{\Fn}(M ,2)} \dot L\in {[ J]}^{{\omega}_1},
      \end{equation}
      and for ${\alpha}\in J$ and $p' \in {\Fn}({\omega}_1,2)$,
      \begin{equation}\label{eq:qlepalpha}
         \text{$p' \Vdash_{{\Fn}(M,2)} {\alpha}\in \dot L$ iff $p' \le p_{\alpha}$.
         }
      \end{equation}

       Since $\{y_{\alpha}:{\alpha}\in {\omega}_1\}\in {[{\omega}_1]}^{{\omega}_1}$ and $1 \Vdash_{{\Fn}({\kappa},2)} \<\omega_1,E\> \in \Gamma({\omega}_1)$,
       there are an  ${\Fn}({\kappa},2)$-name $\{\dot {\mu}_i:i\in {\omega}_1\}$
       and a condition $q \in {\Fn}({\kappa},2)$ with  $q \le p$ such that
      \begin{displaymath}
      q\Vdash_{{\Fn}({\kappa},2)}
      \{\dot {\mu}_i:i\in {\omega}_1\}\in {[\dot L]}^{{\omega}_1},
      \end{displaymath}
      moreover $$q\Vdash_{{\Fn}({\kappa},2)} \forall \{i, j\} \in [\omega_1]^2 \, \{ y_{\dot \mu_i} , {y}_{\dot \mu_j}\} \in E.$$
      For each $i<{\omega}_1$ we may pick $q_i \in {{\Fn}({\kappa},2)}$ with $q_i\le q$ and
      ${\nu}_i \in J$ such that $q_i\Vdash_{{\Fn}({\kappa},2)}\dot {\mu}_i= {\nu}_i$.
      Since $q_i\Vdash_{{\Fn}({\kappa},2)} {\dot \mu}_i = {\nu}_i\in \dot L$, we have $q_i\le p_{{\nu}_i}$ by \eqref{eq:qlepalpha}.
      Thus, we have
      \begin{displaymath}
      q_i\Vdash_{{\Fn}({\kappa},2)} \dot x_{{\dot \mu}_i} =
      \dot x_{{\nu}_i}= y_{{{\nu}_i}}.
      \end{displaymath}
      The following is then our key observation:
      \begin{equation}\label{eq:key}
      \text{if $i \ne j$ and $q_i$ and $q_j$ are compatible, then
      $\{y_{{\nu}_i},y_{{\nu}_j}\} \in E$}.
      \end{equation}
      We may now pick
      $K' \in {[{\kappa}\setm K]}^{{\omega}_1}$ such that
      $q_i\in {\Fn}(K\cup K', 2)$ for all $i<{\omega}_1$.
      Let ${\sigma}:K \cup K' \to {\omega}_1$ be a bijection such that
      ${\sigma}\restriction K = id$.
      Then,  ${\sigma}$ generates a natural order isomorphism ${\sigma}^*$ of ${\Fn}(K \cup K', 2)$
      onto ${\Fn}({\omega}_1,2)$.

      It follows that ${\sigma}^*(q_i) \in Fn ({\omega}_1, 2)$ and
      ${\sigma}^*(q_i)\le {\sigma}^*(q) \le p$ for all $i<{\omega}_1$,
      and $q_i$ and $q_j$ are compatible iff ${\sigma}^*(q_i)$ and ${\sigma}^*(q_j)$ are compatible.
      So, applying \eqref{eq:key} we obtain
      \begin{equation}\label{eq:key2}
         \text{if $i \ne j$ and ${\sigma}^*(q_i)$ and ${\sigma}^*(q_j)$ are compatible, then
         $\{y_{{\nu}_i},y_{{\nu}_j}\}\in E$.}
         \end{equation}


      Now, as is well-known, there is
      $r\in {\Fn}({\omega}_1,2)$ with $r\le {\sigma}^*(q)$  such that
      \begin{equation}\label{eq:r}
      r\Vdash_{{\Fn}({\omega}_1,2)} \dot N = \{i<{\omega}_1:{\sigma}^*(q_i)\in \mc G\} \in {[{\omega}_1]}^{{\omega}_1},
      \end{equation}
      where $\mc G$ is the canonical name of the ${\Fn}({\omega}_1,2)$-generic filter.
      Then, by \eqref{eq:key2}, we have
      \begin{equation}\label{eq:key3}
      \text{if $i \ne j$ and ${\sigma}^*(q_i), {\sigma}^*(q_j) \in \mc G$  then
         $\{y_{{\nu}_i},y_{{\nu}_j}\}\in E$.}
      \end{equation}
      Since ${\sigma}^*(q_i) \le p_{{\nu}_i}$,
      \begin{equation}\label{eq:key4}
      {\sigma}^*(q_i) \Vdash_{{\Fn}({\omega}_1, 2)} i\in \dot N\,\Rightarrow\,(\dot {\mu}_i= {\nu}_i\,\wedge\,
      \dot x_{{\nu}_i} = y_{{\nu}_i} ).
      \end{equation}
      Thus, putting together \eqref{eq:r}, \eqref{eq:key3} and \eqref{eq:key4}, we conclude that
      \begin{displaymath}
         r \Vdash_{{\Fn}({\omega}_1,2)}|\dot N|={\omega}_1 \land
         \big[\{\dot x_{{\nu}_i}:i\in \dot N\}\big]^2\subs E.
         \end{displaymath}
       This completes our proof.
      \end{proof}

\begin{theorem}\label{tm:cohen3}
   If $X$ is any space in our ground model $V$ with $|X| = \weight(X) = \omega_1$ then
   $$V^{{\Fn}({\kappa},2)}\models X\in C({\omega}_1)$$
  implies
   $$V^{{\Fn}({\omega}_1,2)}\models X\in C({\omega}_1).$$
   \end{theorem}

\begin{proof}[Proof of Theorem \ref{tm:cohen3}]
Fix in $V$ a base $\mc B$ of $X$ with $|\mc B|= \omega_1$ and
let $S = \{\<x,B\> : x\in B \in \mc B \}$, moreover
let
\begin{displaymath}
E=\{\{\<x_0,B_0\>,\<x_1,B_1\>\}\in {[S]}^{2}: \{x_0,x_1\}\subs B_0\cap B_1\}.
\end{displaymath}
Clearly, then
$$X \in C({\omega}_1) \,\Leftrightarrow\,\<S, E\> \in \Gamma({\omega}_1).$$
But $\mc B$ remains a base in any generic extension of $V$, hence
Theorem \ref{tm:cohen-graph} applied to the graph $\<S, E\>$ together with the above equivalence
completes our proof.
\end{proof}

\begin{theorem}\label{tm:cohen2}
   If $V \models CH$  then, for every cardinal $\kappa$,
\begin{displaymath}
   V^{{\Fn}({\kappa},2)}\models (X\in C({\omega}_1) \cap T_3\,\wedge\,
   |X|=\weight(X)={\omega}_1) \, \Rightarrow\,  \nw(X)={\omega}.
\end{displaymath}
\end{theorem}

\begin{proof}[Proof of Theorem \ref{tm:cohen2}]
Fix in $V^{{\Fn}({\kappa},2)}$ a base  $\mc B$ of $X$ with $|\mc B|= \omega_1$. Then there is some $I\in V \cap {[{\kappa}]}^{{\omega}_1}$
such that both $\mc B$ and the underlying set of $X$ belong to $V^{{\Fn}(I,2)}$.
We may then apply Theorem \ref{tm:cohen3} for $X$ in $V^{{\Fn}(I,2)}$
to obtain that
\begin{displaymath}
V^{{\Fn}(I,2)*{\Fn}({\omega}_1, 2)}\vDash X\in C({\omega}_1) \cap T_3.
\end{displaymath}
But $CH$ holds in  $V^{{\Fn}(I,2)*{\Fn}({\omega}_1,2)}$ as well, consequently
we have $$V^{{\Fn}(I,2)*{\Fn}({\omega}_1,2)}\models \nw(X)={\omega}$$
by Corollary \ref{co:kík}. But then we have
$V^{{\Fn}({\kappa},2)}\models \nw(X)={\omega}$ as well,
because a network remains a network in any generic extension.
\end{proof}

Now, we are ready to prove  Theorem \ref{tm:cohen}.

\begin{proof}[Proof of Theorem \ref{tm:cohen}]
We are actually going to prove the following stronger result:

\smallskip

{\em Assume that
\begin{displaymath}
   V^{{\Fn}({\kappa},2)} \models  X\in C({\omega}_1) \cap T_3,
   \end{displaymath}
moreover $Y \in [X]^{\omega_1}$. Then for any topology $\varrho$ on $Y$ that is coarser than
its subspace topology, i.e. $\varrho \subs \tau(Y) = \tau(X) \upharpoonright Y$, and has weight $\,\weight(\varrho) \le \omega_1$,
we have $\nw(\varrho) = \omega$.}

This, of course is an immediate consequence of Theorem \ref{tm:cohen2} if $\varrho$ is $T_3$ because
$\varrho \in C(\omega_1)$ is obvious. We next show that the general case can be reduced to the $T_3$ case.

To do this, we first fix a base $\mc B$ of $\varrho$
with $|\mc B| \le \omega_1$.
We then take an elementary submodel $M$
of $H(\vartheta)$ for a sufficiently large regular cardinal $\vartheta$ such that $Y, \mc B \in M$ and $|M| = \omega_1$.
Note that $|Y| = \omega_1 \ge |\mc B|$ imply $Y  \cup \mc B \subs M$.

Now, consider the topology $\tau_M$ on $Y$ generated by $M \cap \tau(Y)$.
Theorem 1.1 of \cite{JunTall98} says that then $\tau_M \in T_3$, hence $\tau_M \in C(\omega_1) \cap T_3$ and $\weight(\tau_M) \le \omega_1$ imply
$\nw(\tau_M) = \omega$. But $\varrho \subs \tau_M$ implies that a network for $\tau_M$ is also a network for $\varrho$,
hence we have $\nw(\varrho) = \omega$, completing the proof.
\end{proof}

\section{Some separation results}

In this section we are going to use forcing to obtain two consistent examples of 0-dimensional $T_2$, hence $T_3$
spaces. The first example separates $N(\omega_1)$ and $K(\omega_1)$, while the second separates $K(\omega_1)$
from $C(\omega_1)$.

We start with a subsection on a technical tool that will play a crucial role in our forthcoming forcing constructions.

\mysubsection{Good colorings}

We first fix some notation.
If $M < \omega$ and $A$ is any $M$-element set of ordinals, we let $\{A(i):i<M\}$ be the increasing enumeration of $A$.
Following the notation from \cite{Juhasz2002HFD}, we shall denote by $\mc D^M_{{\omega}_1}({\eta})$ the collection of all
families $\{A_{\alpha}:{\alpha}<{\omega}_1\}\subs [\eta]^M$, whose members are pairwise disjoint.

\begin{definition}\label{df:good}
Given an uncountable ordinal ${\eta}$, we call a function $g:{\llambda}\times {\eta}\to 2$ {\em good}
iff
\begin{enumerate}
\item
for any $M\in {\omega}$, ${\varepsilon}:M^2\to 2$,
$\{A_{\alpha}:{\alpha}<{\omega}_1\}\in \mc D^M_{{\omega}_1}({\llambda})$ and\\
$\{B_{\alpha}:{\alpha}<{\omega}_1\}\in \mc D^M_{{\omega}_1}({\eta})$,
there are ${\alpha}<{\beta}<{\omega}_1$ such that
\begin{equation}\label{eq:good}
\forall i<M\  \forall j<M\quad
g(A_{\alpha}(i),B_{{\beta}}(j))
=g(A_{{\beta}}(i),B_{\alpha}(j))={\varepsilon}(i,j),
\end{equation}
\item
for any $\{\alpha, \beta\} \in [\omega_1]^2$ there is $\xi < \eta$ with $g(\alpha, \xi) \ne g(\beta, \xi)$.

\end{enumerate}
\end{definition}

\begin{definition}\label{df:spaceeta}
   Assume that
   ${\eta}$  is an  uncountable ordinal, and
   $g:{\llambda}\times {\eta}\to 2$ is a function.
For ${\zeta}<{\eta}$ and $i<2$ write

 \begin{displaymath}
 U_g({\zeta},i)=\{{\alpha}<\llambda: g({\alpha},{\zeta})=i\}.
 \end{displaymath}
Then, we denote by ${\tau}_g$  the topology
 on ${\llambda}$ generated by $\{U_g({\zeta},i):{\zeta}<{\eta},i<2\}$.
This means that all finite intersections of members of $\{U_g({\zeta},i):{\zeta}<{\eta},i<2\}$
form a base $\mc B_g$ for ${\tau}_g$. The members of $\mc B_g$ are then of the form
$$B_g[\varepsilon] = \bigcap \{U_g(\zeta, \varepsilon(\zeta)) : \zeta \in \dom(\varepsilon)\},$$
with $\varepsilon$ ranging over ${\Fn}(\eta, 2)$.

We say that $g$ is {\em independent} iff every non-empty $U\in {\tau}_g$ has cardinality ${\omega}_1$,
i.e. $\Delta(\tau_g) = \omega_1$.
\end{definition}

It is clear that ${\tau}_g$ is 0-dimensional and $T_2$ with $\Delta({\tau}_g)={\omega}_1$ for any independent good function $g$.
Now, the following proposition is straight forward to check.

\begin{proposition}\label{pr:good}
If $g$ is any ${\Fn}(\omega_1 \times \eta, 2)$-generic function then $g$ is both good
and independent in $V^{{\Fn}(\omega_1 \times \eta, 2)}$.
\end{proposition}

Next, given a good function $g$, we define a forcing notion
$\<Q_g,\le \>$ that will be used to force a countable network for ${\tau}_g$. The underlying set
$Q_g$ will be a subset of ${\Fn}({\omega}\times{\llambda},2)\times {[{\omega}\times {\eta}\times 2]}^{<{\omega}}$,
defined as follows.
For  $q=\<s,\mc F\>\in {\Fn}({\omega}\times{\llambda},2)\times {[{\omega}\times {\eta}\times 2]}^{<{\omega}}$ we
 write $$N^q(\ell)=\{{\alpha}\in {\llambda}: s(\ell,{\alpha})=1\}.$$ Then we set
 \begin{multline*}\label{eq:star}
Q_g=\Big\{\<s,\mc F\>\in {\Fn}({\omega}\times{\llambda},2)\times {[{\omega}\times {\eta}\times 2]}^{<{\omega}}:\\
\text{ $\<\ell,{\zeta},i\>\in \mc F$ implies $N^q(\ell)\subs U_g({\zeta},i)$}\Big\}.
\end{multline*}
Clearly, we have $|Q_g| = |\eta|$.

Then we simply put $\<s_1,\mc F_1\>\le \<s_0\,\mc F_0\>$ iff $s_1\supset s_0$ and $\mc F_1\supset \mc F_0$.

For any $q=\<s,\mc F\>\in Q_g$ we put
\begin{displaymath}
   A(q)=\{{\alpha}<\llambda:\exists \ell<{\omega}\ \<\ell,{\alpha}\>\in \dom(s)\},
\end{displaymath}
and
\begin{displaymath}
   B(q)=\{{\zeta}<{\eta}:\exists \ell<{\omega}\ \exists i<2 \<\ell,{\zeta},i\>\in \mc F\}.
\end{displaymath}

\begin{lemma}\label{lm:w-net}
$1_{Q_g}\Vdash \nw(\tau_g)={\omega}$.
\end{lemma}

\begin{proof}
Given a generic filter $\mc G\subs Q_g$,  for $\ell<{\omega}$ let
\begin{displaymath}
N_\ell=\bigcup\{N^q(\ell):q\in \mc G\}.
\end{displaymath}
Now, if $g({\alpha},{\zeta})=i$ (i.e. ${\alpha}\in U_g({\zeta},i)$), then
\begin{displaymath}
D_{{\alpha},{\zeta}, i}=\big\{\<s,\mc F\>\in Q_g:  \exists \ell\in {\omega}\
(s(\ell,{\alpha})=1 \land \<\ell,{\zeta},i\>\in \mc F) \big\}
\end{displaymath}
is clearly dense in $Q_g$.
Hence,
$$V[\mc G]\models \forall {\alpha}<{\omega}_1\ \forall {\zeta}<{\eta}\ \forall i<2\,\big(
{\alpha}\in U_g({\zeta},i)\,\Rightarrow\,  \exists \ell\in {\omega}\ ({\alpha}\in N_\ell\subs U_g({\zeta},i))\big).$$
But clearly, then all the finite intersections of the $N_\ell$'s form a countable network for ${\tau}_g$ in $V[\mc G]$.
\end{proof}

We actually haven't used so far that the function $g:{\llambda}\times {\eta}\to 2$ is good.
This, however will play an essential role in the key lemma below that will provide us
the c.c.c. property  of all finite support powers of $Q_g$.

Before presenting this lemma, we need a preparatory fact about $\Delta$-systems  that is most likely well-known
but we could  not find it in the literature. So we decided to include it here together with a proof.

\begin{fact}\label{fc:stacdelta}
            If
            $\{E_{\nu}:{\nu}\in D\}$ is a family of finite sets indexed by a stationary set $D\subs {\omega}_1$ then
            there is a stationary subset $D'$ of $D$ such that
            $\{E_{\nu}:{\nu}\in D'\}$ forms a $\Delta$-system.
            \end{fact}

\begin{proof}[Proof of the Fact]
We can assume that $\bigcup_{{\nu}\in D}E_{\nu}\subs {\omega}_1$ because $|\bigcup_{{\nu}\in D}E_{\nu}| \le {\omega}_1$.
Next, consider the following club subset $C$ of ${\omega}_1$: $$C=\{{\gamma}<{\omega}_1 : \forall\,{\beta} \in D \cap {\gamma}\,(E_\beta \subs \gamma)\}.$$
Clearly then for any $\{\nu, \mu\} \in [D \cap C]^2$ we have that $E_\nu \setm \nu$ and $E_\mu \setm \mu$ are disjoint.

As $D \cap C$ is stationary, there are a $k < \omega$ and a stationary $D^* \subs D \cap C$ such that
$|\nu \cap E_\nu| = k$ for all $\nu \in D^*$. We may now apply Fodor's pressing down theorem $k$ times
to obtain a stationary $D' \subs D^*$ and a $k$-element subset $E$ of $\omega_1$ such that
$E = \nu \cap E_\nu$ for all $\nu \in D'$. But then $\{E_{\nu}:{\nu}\in D'\}$ is indeed a $\Delta$-system with kernel $E$.
\end{proof}

Now we are ready to present our key lemma promised above.

\begin{lemma}\label{lm:key-lemma}
   Assume that $g : \omega_1 \times \eta \to 2$ is a good function, moreover
   $M < {\omega}$,
   $W\subs {\omega}_1$ is stationary, and
   $$\{q_{\nu} = \<q_{{\nu},m}:m<M\> : {\nu}\in W\} \subs (Q_g)^M$$ is
   such that
   \begin{equation}\label{eq:WAnu}
   W\cap \bigcup\{A(q_{{\nu},m}):{\nu}\in W, m<M\}=\empt.
   \end{equation}
   Then for any $e < 2$ there is a pair  $\{{\nu},{\mu}\}\in {[W]}^{2}$ such that
   $q_{\nu}$ and $q_{\mu}$ are compatible in $(Q_g)^M$ and $g({\nu},{\mu})=e$.
   \end{lemma}

   \begin{proof} Write $q_{{\nu},m}=\<s_{{\nu},m},\mc F_{{\nu},m}\>$
      for ${\nu}\in W$ and $m<M.$
      For ${\nu}\in W$   put
         $A_{{\nu}}=\bigcup\{A(q_{{\nu},m}):m<M\}$, $B_{\nu}=\bigcup\{B(q_{{\nu},m}):m<M\}$
         and
      $C_{\nu}=\{{\nu}\}\cup A_{\nu}\cup B_{\nu}$.

         By appropriate thinning out and using Fact \ref{fc:stacdelta} three times, we can assume that
         \begin{enumerate}[(1)]
            \item The families $\{A_{\nu}:{\nu}\in W\}$,
         $\{B_{\nu}:{\nu}\in W\}$, and $\{C_{\nu}:{\nu}\in W\}$ form $\Delta$-systems with
         kernels $A$ and $B$ and $C$ respectively.
         \item For each $\{{\nu},{\mu}\}\in {[W]}^{2}$ we have $|C_{\nu}|=|C_{\mu}| = K$ and the  order preserving bijection
         ${\sigma}_{{\nu},{\mu}}:C_{\nu}\to  C_{\mu}$ has the following properties:
         \begin{enumerate}[(i)]
         \item ${\sigma}_{{\nu},{\mu}}\restriction C$ is the identity, \smallskip
         \item ${\sigma}_{{\nu},{\mu}}({\nu})={\mu}$, \smallskip
         \item $s_{\mu,m}=\{\<\<\ell,{\sigma}_{{\nu},{\mu}}({\alpha})\>,i\>:\<\<\ell,{\alpha}\>,i\>\in s_{{\nu},m}\}$ for each $m<M$,\smallskip
         \item $\mc F_{\mu,m}=\{\<\ell,{\sigma}_{{\nu},{\mu}}({\zeta}),i\>:\<\ell,{\zeta},i\>\in \mc F_{\nu,m}\}$ for each $m<M$,\smallskip
         \item $g({\alpha},{\beta})=g({\sigma}_{{\nu},{\mu}}({\alpha}),{\sigma}_{{\nu},{\mu}}({\beta}))$ for each $\<{\alpha},{\beta}\>\in (C_{\nu})^2\cap \dom(g)$.\smallskip
         \end{enumerate}
         \end{enumerate}

      We next  fix a function  ${\varepsilon}:K\times K\to 2$ such that:
      \begin{enumerate}[(a)]
      \item if ${\nu}=C_{\nu}(k)$, and hence ${\mu}=C_{\mu}(k)$, then ${\varepsilon}(k,k)=e$,
      \item for each $i,j<K$, if $C_{\nu}(i)\in A_{\nu}$ and $C_{\nu}(j)\in B_{\nu}$, then
      \begin{equation}\label{eq:epsilong}
         {\varepsilon}(i,j)=g(C_{\nu}(i),C_{\nu}(j)).
      \end{equation}
      \end{enumerate}
      Since ${\nu}\notin A_{\nu}$ by our assumption, we can find such an ${\varepsilon}$.
         The function  $g$ is good, so there is a pair  $\{{\nu},{\mu}\}\in [W]^2$ such that
         \begin{equation}\label{eq:good-appl}
            \forall i<K\  \forall j<K\quad
            g(C_{\nu}(i),C_{{\mu}}(j))
            =g(C_{{\mu}}(i),C_{{\nu}}(j))={\varepsilon}(i,j).
            \end{equation}
      Thus, (a) and  \eqref{eq:good-appl} imply
      \begin{displaymath}
      g({\nu},{\mu})=e,
      \end{displaymath}
      and (b) and  \eqref{eq:good-appl} imply
         \begin{equation}\label{eq:good-2}
         \forall {\alpha}\in A_{\nu}\ \forall {\xi}\in B_{\nu}
         \ g({\alpha},{\sigma}_{{\nu},{\mu}}({\xi}))=g({\sigma}_{{\nu},{\mu}}({\alpha}),{\xi})=g({\alpha},{\xi}).
         \end{equation}
         We claim that $r_m = \<s_{\nu,m}\cup s_{\mu,m},\mc F_{\nu,m}\cup \mc F_{\mu,m}\>\in Q_g$ for each $m<M$
         and so $\<r_m : m < M\>$ is a common extension of $q_{\nu}$ and $q_{\mu}$.

         Since  $s_{\nu,m}\cup s_{\mu,m}\in {\Fn}({\omega}\times{\eta},2)$ is obvious,
we should only check that
      \begin{displaymath}
         \text{$\<\ell,{\zeta},i\>\in \mc F^{r_m}$ implies $N^{r_m}(\ell)\subs U_({\zeta},i)$}.
      \end{displaymath}
So, assume that ${\alpha}\in N^{q_m}(\ell)$, i.e.
$(s_{\nu,m}\cup s_{\mu,m})(\ell,{\alpha})=1$.
By symmetry, we can assume that ${\alpha}\in A(q_{{\nu},m})$.
If ${\zeta}\in B(q_{{\nu},m})$ then we are done because $q_{{\nu},m} \in Q_g$.
So we can assume that $\<\ell,{\zeta},i\>\in \mc F_{\mu,m}$.

         Let us put ${\zeta}^*={\sigma}_{{\nu},{\mu}}^{-1}({\zeta})$.
         Then, by (iii), $\<\ell,{\zeta}^*,i\>\in \mc F_{{\nu},m}$, hence
         $s_{\nu,m}(\ell,{\alpha})=1$
         implies  $g({\alpha},{\zeta}^*)=i$.
         Then, applying \eqref{eq:good-2} we obtain $g({\alpha},{\zeta})=g({\alpha},{\zeta}^*) = i$,
         i.e. ${\alpha}\in U_g({\zeta},i)$.
This completes the proof.
      \end{proof}

\begin{lemma}\label{lm:ccc}
If $g:{\omega}_1\times {\eta}\to 2$ is a good function, $M\in {\omega}$, then
the poset $(Q_g)^M$ satisfies c.c.c. Consequently, any finite support power of $Q_g$ satisfies c.c.c.
\end{lemma}
\begin{proof}
Assume that $$\{q_{\nu} = \<q_{{\nu},m}:m<M\> : {\nu}\in \omega_1\} \subs (Q_g)^M,$$
where $q_{{\nu},m}=\<s_{{\nu},m},\mc F_{{\nu},m}\>$.
We may pick $I\in {[{\omega}_1]}^{{\omega}_1}$ such that
$$A=\bigcup\{A(q_{{\alpha},m}):{\alpha}\in I,m<M\}\in NS({\omega}_1).$$ Let $W={\omega}_1\setm A$ and
reindex $I$ in a one-one way by $W$ as $I=\{{\alpha}_{\nu}:{\nu}\in W\}$.

We can then apply  Lemma \ref{lm:key-lemma} for
$\<q_{{\alpha}_{\nu}}:{\nu}\in W\>$
to find $\{{\nu},{\mu}\}\in {[W]}^{2}$
such that $\alpha_\nu \ne \alpha_\mu$ and $q_{\alpha_{\nu}}$, $q_{{\alpha}_{\mu}}$ are compatible.
Now, the second part follows immediately.
\end{proof}

\mysubsection{Forcing a space in $K(\omega_1) \setm N(\omega_1)$}

First we present two simple and  easy preparatory lemmas.

\begin{lemma}\label{lm:nww1}
Let $\tau$ be a topology on $\omega_1$ and assume that for every countable family
$\mc A \subs [\omega_1]^{\omega_1}$ there is $U \in \tau$ with $|U| = \omega_1$
such that $A \not\subs U$ for all $A \in \mc A$. Then $\nw(\tau) = \omega_1$.
\end{lemma}

\begin{proof}
Indeed, if $\mc N$ is any countable family of subsets of $\omega_1$
then there is $\alpha < \omega_1$ such that $A \subs \alpha$ for every countable $A \in \mc N$.
Then, by our assumption,
there is an uncountable $U \in \tau$ such that $A \not\subs U$ whenever $A \in \mc N$ is uncountable.
Now, if $x \in U \setm \alpha$
and $x \in A \in \mc N$ then $A \not\subs U$, hence $\mc N$ is not a network for $\tau$.
\end{proof}

\begin{lemma}\label{lm:splitting}
   Assume that
   $g:{\omega}_1\times {\lambda}\to 2 $ is a  good function.
   Then for each $A \in  {[{\omega}_1]}^{{\omega}_1}$ the set
\begin{displaymath}
 K(A)=\{{\zeta}<{\lambda}: A \subs U_g({\zeta},0)\}
\end{displaymath}
is countable.
\end{lemma}

\begin{proof}
Consider any set
$B\in {[{\lambda}]}^{{\omega}_1}$ of size $\omega_1$.
Since $g$ is good, we may apply item (1) of Definition \ref{df:good} (with $M = 1$)
for the sets of singletons $\{\{a\} : a \in A\}$ and $\{\{b\} : b \in B\}$
to obtain $a \in A$ and $b \in B$
such that $g(a, b) = 1$, hence $a \in A \setm U_g(b, 0)$.
Since $B$ was an arbitrary subset of ${\lambda}$ of cardinality $\omega_1$,
it follows that $K(A)$ is indeed countable.
\end{proof}

\begin{theorem}\label{tm:main}
   Assume that ${\lambda}$ is a  cardinal with  $cf({\lambda})>{\omega}_1$ and
   $g:{\omega}_1\times {\lambda}\to 2 $ is an independent good function.
Then the finite support product poset  $$Q = \prod_{\omega_1 \le {\eta}<{\lambda}}Q_{g\restriction {\omega}_1\times {\eta}}$$
is c.c.c, moreover we have
   $$V^Q\models  X_g = \<\omega_1, \tau_g\> \in K({\omega}_1)\setm N({\omega}_1).$$
   \end{theorem}

\begin{proof}[Proof of Theorem \ref{tm:main}]
$Q$ is clearly a regular subposet of
the finite support power $Q_{g}^{\lambda}$ of $Q_g$, hence it is c.c.c.
by Lemma \ref{lm:ccc}.

To see that $V^Q\models  X_g \in K({\omega}_1)$,
consider in $V^Q$ any $U \in PNA(X_g, \omega_1)$  such that
$$U(\alpha) = B_g[\varepsilon_\alpha] \in \mc{B}_g$$
for all $\alpha \in \dom(U)$, with $\varepsilon_\alpha \in {\Fn}(\lambda, 2)$. Now, $cf({\lambda})>{\omega}_1$ implies that
there is ${\eta}<{\lambda}$ such that $U \in PNA(X_{g\restriction {\omega}_1\times {\eta}}, \omega_1)$.
Then $V^Q\models \nw({g\restriction {\omega}_1\times {\eta}})={\omega}$ by Lemma \ref{lm:w-net}, so we have
$\nw(X_U)={\omega}$ in $V^Q$ as well. Thus we indeed have $X_g\in K({\omega}_1)$ in $V^Q$.

Next we shall show that $V^Q\models  \nw(X_g) = \omega_1$, hence $V^Q\models  X_g \notin N(\omega_1)$.
To do this, we first fix a $Q$-generic filter $\mc G\subs Q$ over $V$.
In view of Lemma \ref{lm:nww1}, $V^Q\models  \nw(X_g) = \omega_1$ immediately follows from the following claim.

\begin{claim}\label{cl:w1}
For any countable family $\mc A = \{A_n : n \in \omega\} \subs [\omega_1]^{\omega_1}$ in $V[\mc G]$, there is $\eta < \lambda$ such that $A_n \not\subs U_g(\eta, 0)$ for all $A_n \in \mc A$.
\end{claim}

Now, let us turn to the proof of Claim \ref{cl:w1}.
Since $cf({\lambda})>{\omega}_1$, there is ${\mu}<{\lambda}$
such that $\mc A\in V[\mc G']$,
where $Q' = \{q \in Q : \supp(q) \subs \mu\}$ and $\mc G'=\mc G \cap Q'$.

For any $n<{\omega}$ let $\dot A_n$ be a $Q'$-name for $A_n$, and for notational simplicity
assume that $${1}_{Q'} \Vdash' \{\dot A_n:n\in {\omega}\} \subs [\omega_1]^{\omega_1}.$$
Of course, $\Vdash'$ denotes forcing with $Q'$.

Now, for any $n < \omega$ and $p\in Q'$ set
\begin{displaymath}
B(n,p)=\{{\zeta}<{\omega}_1: \exists\ q\le p\,\, q \Vdash' {\zeta}\in \dot A_n\},
\end{displaymath}
and put
\begin{displaymath}
\mc B=\{B(n,p) : n<{\omega},\, p\in Q'\}.
\end{displaymath}
Note that we have $\mc B \in V$, moreover
$p\Vdash' \dot A_n\subs B(n,p)$ implies that
$\mc B\subs {[{\omega}_1]}^{{\omega}_1}$. Let
$K=\bigcup\{K(B): B\in  \mc B\}$, then
by lemma \ref{lm:splitting} we have
$|K|\le |\mc B|\cdot {\omega}\le |Q'|\cdot {\omega}<{\lambda}$,
so we can pick ${\eta}\in {\lambda}\setm K$, and then
\begin{equation}\label{eq:delta}
\forall n<{\omega}\ \forall p\in Q'\
B(n,p)\not\subs U_g({\eta},0).
\end{equation}
So, by \eqref{eq:delta}, for every $n \in \omega$ and
$p\in Q'$
there are $q\le p$ and ${\zeta}<{\omega}_1$
such that $q\Vdash' {\zeta}\in \dot {A}_n \setm U_g(\eta,0)$.
Hence the set
\begin{displaymath}
E_n=\{q\in Q':q'\Vdash \dot A_n \not\subs U_g({\eta},0)\}
\end{displaymath}
is dense in $Q'$. Thus $\mc G'\cap E_n\ne \empt$, and so
$V[\mc G']\models A_n \not\subs U_g({\eta},0)$,
completing the proof of the claim and hence of the theorem.
\end{proof}

\mysubsection{Forcing a space in $C(\omega_1) \setm K(\omega_1)$}

We start by fixing in our ground model $V$ an independent good function $g:{\omega}_1\times {\omega}_1\to 2$.
Then, for any $S \subs \omega_1$  we let $$Q(S)=\{q\in Q_g: A(q)\subs S\}=Q_{g\restriction S\times {\omega}_1}.$$
We next consider the finite support product poset
\begin{displaymath}
Q = {\prod} \{Q(S) \,:\, S\in NS({\omega}_1)\}.
\end{displaymath}
Theorem \ref{tm:stac} will then easily follow from our next result.

\begin{theorem}\label{tm:stac-good}
\begin{displaymath}
V^Q\models \{Y\subs {\omega}_1: \nw({\tau}_g\restriction Y)={\omega}\}=NS({\omega}_1).
\end{displaymath}
\end{theorem}

\begin{proof}[Proof of Theorem \ref{tm:stac-good}]
Since $Q$ is again a regular subposet of a finite support power of $Q_{g}$, it satisfies c.c.c.
This implies that if $T\in NS({\omega}_1)\cap V^{Q}$
then there is $S\in NS({\omega}_1)\cap V$ with
$T\subs S$. So, as
\begin{displaymath}
1_{Q(S)}\Vdash
\nw({\tau}_g\restriction S)={\omega}
\end{displaymath}
by Lemma \ref{lm:w-net}, we have $\nw({\tau}_g\restriction T)={\omega}$ in  $V^Q$.

The burden of the proof is to verify that $\nw({\tau}_g\restriction Y)={\omega}_1$
whenever $Y\subs {\omega}_1$ is stationary in $V^Q$.
Assume, on the contrary, that
\begin{displaymath}
q \Vdash_Q \dot Y\subs {\omega}_1\text{\em \ is stationary and } \{\dot N_\ell:\ell<{\omega}\}
\text{\em\  is a network for }{\tau}_g\restriction \dot Y.
\end{displaymath}
Then the ground model set
\begin{displaymath}
   Z=\{{\nu}<{\omega}_1:\exists q_{\nu}\le q\,  \exists  \<\ell_{\nu}, e_{\nu}\> \in \omega \times 2\ \big(
   q_{\nu}\Vdash_Q  {\nu}\in \dot Y\land {\nu}\in \dot N_{\ell_{\nu}}\subs U_g({\nu},e_{\nu})\big)\}.
   \end{displaymath}
is stationary in ${\omega}_1$ in $V$.
We may then pick $\ell<{\omega}$ and $e<2$ such that
\begin{displaymath}
Z' = \{{\nu}\in Z: \ell_{\nu}=\ell\land e_{\nu}=e\}
\end{displaymath}
is stationary in ${\omega}_1$ as well.
Thus, for each ${\nu}\in Z'$
\begin{displaymath}
q_{\nu}\Vdash_Q {\nu}\in \dot Y\land  {\nu}\in \dot N_\ell \subs U_g({\nu},e).
\end{displaymath}
Consequently, for any $\{{\nu},{\mu}\}\in {[Z']}^{2}$,
\begin{equation}\label{eq:ck}
\text{ if $q_{\nu}$ and $q_{\mu}$
are compatible in $Q$
then $g({\nu},{\mu})=e$},
\end{equation}
because if $q'$ is a common extension of $q_{\nu}$ and $q_{\mu}$,
then
\begin{displaymath}
q'\Vdash_Q {\nu}\in \dot N_\ell \subs U_g({\mu},e),
\end{displaymath}
and so ${\nu}\in U_g({\mu},e)$, i.e. $g({\nu},{\mu})=e$.
Next, however, we'll show that there is a pair $\{{\alpha},{\beta}\}\in {[Z']}^{2}$ such that
$q_{\alpha}$ and $q_{\beta}$ are compatible but $g({\alpha},{\beta})= 1-e$,
giving us the required contradiction.

Since $\dom(q_{\nu})$ is a finite subset of $NS({\omega}_1)$ for all ${\nu} \in Z'$,
we may apply Fact \ref{fc:stacdelta},
to find a stationary set $W_0\subs Z'$ such that
$\{\dom(q_{\nu}):{\nu}\in W_0\}$  forms a $\Delta$-system with kernel $C = \{S_i : i < M\} \subs NS(\omega_1)$.
Thus, for any ${\nu},{\mu}\in W_0$,  the conditions  $q_{\nu}$ and $q_{\mu}$ are compatible iff
$q_{\nu}\restriction C$ and
$q_{{\mu}}\restriction C$ are.
Now, we can apply  Lemma \ref{lm:key-lemma}
for $W = W_0 \setm \cup \{S_i : i < M\}$ and $\{q_{\nu}\restriction C : {\nu}\in W\}$
to find $\{{\nu},{\mu}\}\in {[W]}^{2}$
such that $g({\nu},{\mu})=1-e$ and $q_{\nu}\restriction C$,
$q_{{\mu}}\restriction C$ are compatible.
But then $q_{\nu}$ and $q_{\mu}$ are compatible as well,
contradicting  \eqref{eq:ck}.
\end{proof}

Now, we claim that, in $V^Q$, the space $X_g = \<\omega_1, \tau_g\> \in C(\omega_1) \setm K(\omega_1)$.
Since every set $Y \in [\omega_1]^{\omega_1}$ has a non-stationary subset $S$ with $|S| = \omega_1$,
and so with $\nw({\tau}_g\restriction S)={\omega}$ it is obvious that $X_g \in C(\omega_1)$.
To see that $X_g \notin K(\omega_1)$, we first present the following simple result.

\begin{theorem}\label{tm:U-gen}
If $X$ is any $T_2$-space with $|X|\ge \weight(X)$, then there is
neighborhood assignment $U$ on $X$ such that ${\tau}_U={\tau}(X)$.
\end{theorem}

\begin{proof}
Let $\{B_i:i<{\kappa} = \weight(X)\}$ be a base of $X$.
Pick points $\{y_i, \,z_i:i<{\kappa}\}$ from $X$ that are all distinct. For each $i < \kappa$
let $V_i$ and $W_i$ be disjoint open neighborhoods of $y_i$ and $z_i$, respectively.

Now we define the neighborhood assignment $U$ on $X$ as follows:
\begin{displaymath}
U(x)={}\left\{\begin{array}{ll}
{B_i\cup V_i}&\text{if $x=y_i$},\\
{B_i\cup W_i}&\text{if $x=z_i$},\\
{X} &\text{if $x \in X \setm \{y_i, z_i : i < \kappa\}$.} \\
\end{array}\right.
\end{displaymath}
Then $B_i=(B_i\cup V_i)\cap (B_i\cup W_i)\in {\tau}_U$ for all $i < \kappa$,
hence indeed ${\tau}_U={\tau}(X)$.
\end{proof}
So, by  $|X_g| = \weight(X_g) = \omega_1$ there is $U \in PNA(X_g, \omega_1)$ such that ${\tau}_U={\tau}(X_g)$
hence $\nw(\tau_g) = \omega_1$ implies $X_g \notin K(\omega_1)$.

\medskip

{\em Acknowledgement.} We are grateful to the referee for carefully reading an
earlier version of our paper and, in particular,
spotting a serious error in it.

\end{document}